\documentclass[a4paper]{article}
\usepackage{graphicx}
\vspace{25mm}
\title{Bit Reversal through Direct Fourier Permutation Method and Vectorial Digit Reversal Generalization}

\author{Nicolaie Popescu-Bodorin\footnote{Nicolaie Popescu-Bodorin, Teaching Assistant, Faculty of Mathematics and Computer Science, Spiru Haret University; PhD Student - Computer Science, University of Pitesti; Member of IEEE, IEEE Control Systems Society, IEEE Signal Processing Society; ROMANIA; Email: {\bf bodorin [a] ieee . org}; Correspondence address: {\bf Nicolaie POPESCU-BODORIN, O.P. 19, C.P. 77, Sect.3, Bucharest, ROMANIA.} }}
\date{}
\begin{document}
\maketitle
\pagenumbering{arabic}
\setcounter{page}{1}
% Enter chapter title between curly braces
% Enter section title between curly braces
% Enter subsection title between curly braces
\begin{center}
Proceedings of the 15th 

Conference on Applied and Industrial Mathematics (CAIM), 

 Romanian Society of Applied and Industrial Mathematics (ROMAI), 

12-14 October 2007, Mioveni (Arges), ROMANIA. 
\end{center}
\vspace{20mm}
\begin{abstract}
This paper describes the {\it Direct Fourier Permuation Algorithm}, an  efficient method  of computing Bit Reversal of natural indices $[1,2,3, \dots ,2^{k}]$ in a  vectorial manner ($k$ iterations) and also proposes the {\it Vectorial Digit Reversal Algorithm}, a natural generalization of {\it Direct Fourier Permutation Algorithm} that is enabled to compute the $r$-digit reversal of natural indices $[1,2,3, \dots ,r^{k}]$ where $r$ is an arbitrary radix. Matlab functions implementing these two algorithms and various test and comparative results are presented in this paper to support the idea of  inclusion of these two algorithms in the next Matlab Signal Processing Toolbox official distribution package as much faster alternatives to current Matlab functions $bitrevorder$ and $digitrevorder$. 
\vspace{3mm}


{\it ACM CCS 1998}: G.4., D.3.3.\\

\end{abstract}

\newpage
%\hspace{25mm}
\begin{center}
\section{Introduction}
\end{center}
 \hspace{5mm} Recurrence relations in the Danielson-Lanczos Lemma allow for the immediate implementation of an explicitly recursive function for FFT computation. Each time it calls itself, a  $2^k$ point FFT computation is reduced to two  $2^{k-1}$ FFT computations. This direct approach (divide et impera and explicit backward recursion) allows for the FFT computation algorithm to be expressed by decimation-in-time implementation (1). Subdividing the initial vector $x$ up to a set of pairs on which the FFT computation is very simple is sometimes called {\it buterfly operation}, a name suggested by the graphical representation of the computation [1]. \\

function $X = mkFFT(x)$

\% $x$ - 1x$2^k$ line vector 

\% $X$ - Discrete Fourier Transform of $x$ calculated through 2-radix Decimation  

\%      in Time Fast Fourier Transform with Explicit Backward Recursion.

$N=max(size(x));$

$if N==1$

\hspace{5mm}$X=x;$

else \hspace{110mm} (1)

\hspace{5mm}$oddind = 1:2:N; xodd = x(oddind);$

\hspace{5mm}$evenind = 2:2:N; xeven = x(evenind);$

\hspace{5mm}$O = mkFFT (xodd);$

\hspace{5mm}$ E = mkFFT (xeven);$

\hspace{5mm}For $k=1:1:N/2$

\hspace{5mm}\hspace{5mm}$X(k)=O(k)+(exp(-2*pi*i*(k-1)/N))*E(k);$

\hspace{5mm}\hspace{5mm}$X(N/2+k)=O(k)-(exp(-2*pi*i*(k-1)/N))*E(k);$

\hspace{5mm}end

end;\\

Implicitly, during this computation, a permutation of the argument vector $x$ takes place whenever the function calls itself. Consequently, a composed permutation of the initial vector argument $x$ will occur up to and within the innermost call. This permutation will be referred further in this paper as {\it Fourier Permutation} and is sometimes named the {\it Buterfly Permutation} [1], or {\it Bit Reversal Permutation} [2]-[9] due to the fact that the reversed bit representation of the permuted index of the initial 0-based index (see the following example) is implicitly sorted in ascending order.\\

{\bf Example 1: }
\\
\begin{center}
\begin{tabular}{cccc}

 permuted index &  
 bit reprezentation &
reversed reprezentation &
 0-based index \\ 

$ 0 $ & 
$ 000 $ & 
$ 000 $ & 
$ 0 $ \\

$ 4 $ & 
$ 100 $ & 
$ 001 $ & 
$ 1 $ \\

$ 2 $ & 
$ 010 $ & 
$ 010 $ & 
$ 2 $\\ 

$ 6 $ & 
$ 110 $ & 
$ 011 $ & 
$ 3 $ \\

$ 1 $ & 
$ 001 $ & 
$ 100 $ & 
$ 4 $ \\

$ 5 $ & 
$ 101 $ & 
$ 101 $ & 
$ 5 $ \\ 

$ 3 $ & 
$ 011 $ & 
$ 110 $ & 
$ 6 $ \\

$ 7 $ & 
$ 111 $ & 
$ 111 $ & 
$ 7 $\\

\end{tabular}
\end{center}

\vspace{5mm}

The {\it Fourier Permutation} of the 1-based index $[1,2,3,4,5,6,7,8]$  is $[1,5,3,7,2,6,4,8]$. The stages leading to this permutation through the explicit calls in (1) are the following: $[1,2,3,4,5,6,7,8] \rightarrow [1,3,5,7,2,4,6,8] \rightarrow [1,5,3,7,2,6,4,8] \rightarrow [1,5,3,7,2,6,4,8]$, where a left to right reading shows the permutations operated from the first to the last (the innermost) call.\\

\begin{center}
\section{Additive Constants Method}
\end{center}
 
 \hspace{5mm} In what follows we aim to compute the {\it Fourier Permutation} in an iterative manner, by other methods than the 'Bit Reversal' algorithms proposed in [1] - [9]. 

Our primary intention is to formulate an algorithm for computing the function  $Y_{N} = mkFPerm(N)$, where  $Y_{N}$ is the {\it Fourier Permutation} of the vector of indices $[1,2, \dots ,N]$  and  $N=2^{k}$, with acceptable efficiency. The formulation of such an algorithm requires that a recurrence relation of first order should be found between the consecutive components of the resulting vector $Y_{N} $:

$$
\forall p \in \overline{1,N-1}: Y_{N}(p+1) = Y_{N}(p)+C_{N}(p);\, Y_{N}(1) = 1; \eqno(2)
$$
where $C_{N}(1:N)$  is an additive constants vector to be determined. \\

{\bf Example 2: }

\begin{tabular}{cccccccccc}

$C_{8}(1:7)=[$ & 
$ \, $ & 
$(+4)$ & 
$(-2)$ & 
$(+4)$ & 
$(-5)$ & 
$(+4)$ & 
$(-2)$ & 
$(+4)$ & 
$ ] $ \\ 

$Y_{8}(1:8)=[$ & 
$ 1 $ & 
$ 5 $ & 
$ 3 $ & 
$ 7 $ & 
$ 2 $ & 
$ 6 $ & 
$ 4 $ & 
$ 8 $ & 
$ ] $ \\
\end{tabular}
\\

The {\it Fourier Permutation} of indices $\overline{1,8}$  can be computed taking into account that each component of the resulting vector is the sum of the preceding component and a constant that depends on $N$, and that the first component always has the value 1.

The quantity to be found in all the odd rank positions of vector $C_{N}$  (+4 in the example under discussion) will be further referred to in this paper as the {\it trivial additive constant of the Fourier Permutation}, all the other being called nontrivial constants. \\

{\it Definition} 1. Let $N = 2^{k}, k\in N^{*}$. We shall call the aditive constants of the {\it Fourier Permutation} of indices $[1,2,3, ... , 2^{k}]$ the $(2^{k}-1)$  components of vector $C_{N}$ thus constructed:
\renewcommand{\labelenumi}{\roman{enumi}.}
\begin{enumerate}
\item $C_{2^{1}} = 1$; 
\item $\forall k \in N^{*},\, k \geq 2: \, C_{2^{k}} = [2C_{2^{k-1}}, -2^{k}+3, 2C_{2^{k-1}}]$;
\end{enumerate}

The basic properties of the additive constants corresponding to the {\it Fourier Permutation} of indices $[1,2,3, ... , 2^{k}]$ are given in the following proposition:\\

PROPOSITION 1. 
If $k \geq 2$  and $N = 2^{k}$ , the vector of the additive constants of the {\it Fourier Permutation} ($C_{N}$ ) has the following properties:
\renewcommand{\labelenumi}{\roman{enumi}.}
\begin{enumerate}
\item All odd-rank components store the value of the trivial additive constant: 

 $\forall p \in \overline{1,N/2}: C_{N}(2p-1)=2^{k-1}$;
\item All non-trivial additive constants are negative: 

 $\forall p \in \overline{1,N/2}: C_{N}(2p)<0 $;
\item The lowest non-trivial additive constant splits the $ C_{N}$ vector into two equal vectors: 

 $ C_{N}(1 : 2^{k-1}-1) = C_{N}(2^{k-1}+1 : 2^{k}-1)$;
\item The value of the minor nontrivial constant: 

 $ C_{N}( 2^{k-1}) = -2^{k}+3$;
\item All non-trivial additive constants, except the minor nontrivial constant, are even: 

 $\forall p \in \overline{1,(N/2-1)}: |C_{N}(2p)|\, mod\, 2 = 0$;
\item The components that are symetrically placed in relation to the minor nontrivial constant are equal: 

 $\forall p \in \overline{1,(N/2-1)}: C_{N}(2^{k}-p) = C_{N}(p)$; 
\item Forward recursion: at step $(k+1)$  each of the vectors formed with the components on the left and on the right of the minor non-trivial additive constant, respectively, is double the vector of constants computed at step k: 

 if $N_{1} = 2^{k}$  and $N_{2} = 2^{k+1}$  then: 
 $C_{N_{2}}(1 : N_{1}-1) = 2 C_{N_{1}} = C_{N_{2}}( N_{1}+1 : N_{2})$;
\item Forward recursion between minor nontrivial constants: 

 if $N_{1} = 2^{k}$  and $N_{2} = 2^{k+1}$  then: 
 $C_{N_{2}}(N_{1}) = 2 C_{N_{1}}(2^{k-1})$; 
\end{enumerate}

\vspace{5mm}

{\it Consequence: } 
If $N = 2^{k}$  and $k \geq 2$ , then the number of unique non-trivial additive constants of the {\it Fourier Permutation} is  $(k-1)$, and the value of the trivial additive constant is  $2^{k-1}$. Consequently the computation of the {\it Fourier Permutation} of indices  $[1,2,3, ... , 2^{k}]$ is reduced to:
\begin{itemize}
\item the computation of these k additive constants and their distribution in a template vector of length $2^{k}-1$.
\item the computation of each component of the resulting vector as the sum of the preceding component and the corresponding additive constant.
\end{itemize}

\begin{center}
\section{Computing Fourier Permutation through Additive Constants Algorithm}
\end{center}

\hspace{5mm}Due to property (P1.iii), in order to obtain the template vector of additive constants it suffices to determine its first $2^{k-1}$  components, i.e. its first $2^{k-2}$  non-trivial additive constants (as all the others, i.e. all odd rank components store the trivial additive constant). According to the above considerations, the matlab function for generating the additive constants of the Fourier permutation and the function for computing the {\it Fourier Permutation} of indices $[1,2,3, ... , 2^{k}]$  can be written as: \\

function $V=mkCAPF(N)$

\% $VI =$ intermediate vector of the first $2^{k-1}$ non-trivial additive constants 

\% $N = 2^{k}, k>2$

\% $ct =$ trivial constant

$ct=N/2$; \% property (P1.i)

$k=log2(N)$;

$VI=[ -2 , -5 ]$; 

\% the two non-trivial additive constants corresponding to case k=3

for $p=4:k$

\hspace{5mm}$VI=2*VI$; \% property (P1.vii)

\hspace{5mm}$c=max(size(VI))$;

\hspace{5mm}$VI=[ VI , VI(1:c-1)]$; \% property (P1.iii)

\hspace{5mm}$VI=[ VI , VI(c)-3 ]$; \% properties (P1.vii, P1.viii)

end;

$c=max(size(VI))$;

$VI=[ VI , VI(1:c-1) ]$; \% property (P1.iii)

$V=zeros(1,N-1)+ct$; \% property (P1.i)

$V(2:2:N-1)=VI$; \\

function $V=mkFPerm(N)$ 

$V=zeros(1,N); V(1)=[1];$

$CN=mkCAPF(N)$;

for $p=1:N-1$

\hspace{5mm}$V(p+1) = V(p)+CN(p)$;

end;\\

The complexity of the computation of permuted indices depends on the multiplication operations in the mkCAPF function and on the addition iterated in the mkFPerm function. 
Consequently the complexity of computing function mkFPerm is $O(N log_{2}(N))$  and may decrease if multiplications are excluded from the computational mechanism and the number of iterations decreases, possibly to $k$. \\

{\centering
\section{Direct Fourier Permutation Method}
}

PROPOSITION 2. If  $N = 2^{k}$ and  $k \geq 2$, the vector of the additive constants of the {\it Fourier Permutation}, $C_{N} = mkCAPF(N)$, has the following properties:
\renewcommand{\labelenumi}{\roman{enumi}.}
\begin{enumerate}
\item The sum of all the additive constants of the {\it Fourier Permutation} is: 
$$\sum_{i=1}^{N-1}C_{N}(i) = 2^{k}-1;$$
\item The sum of all non-trivial constants of the {\it Fourier Permutation} is: 
$$\sum_{i=1}^{N/2}C_{N}(2i-1) = -( 2^{k-1}-1)^{2};$$
\item The sum of the first $2^{k-1}$  constants of the {\it Fourier Permutation} is 1: 
$$\sum_{i=1}^{N/2}C_{N}(i) = 1;$$
\item The sum of any $2^{k-1}$ consecutive constants of the {\it Fourier Permutation} is 1: 
$$\forall p \in \overline{1,N/2}: \sum_{i=p}^{N/2+p-1}C_{N}(i) = 1;$$
\end{enumerate}

\vspace{5mm}

{\it Consequences:}

By construction, for any $i \in \overline{1,2^{k-1}}$  the difference between the rank  $(i+2^{k-1})$ component and the rank $i$ component of vector $Y$ is the sum of the  $2^{k-1}$ consecutive additive constants starting with (and including) the rank $i$ constant. \\
Consequently:  
$$Y_{N}(2^{k-1}+1 : 2^{k}) = Y_{N}(1 : 2^{k-1})+1;$$

Moreover, using properties (P1.vii)  and (P2.iii,iv)  applied to case  $N = 2^{k-1}$, it follows that the sum of any  $2^{k-2}$ consecutive constants selected from the first  $(2^{k-1}-1)$ components of  $C_{N}$ is 2. However, by construction, for any $i \in \overline{1,2^{k-2}}$ , the difference between the rank $(i + 2^{k-2})$  component and the rank $i$ component of vector $Y_{N}$ is the sum of the $2^{k-2}$  consecutive additive constants starting with the rank $i$ constant:
$$Y_{N}(2^{k-2}+1 : 2^{k-1}) = Y_{N}(1 : 2^{k-2})+2;$$
 
The sum of any $2^{k-3}$  consecutive constants selected from the first $(2^{k-2}-1)$  components of  $C_{N}$ is $2^{2}$ , and therefore:
$$Y_{N}(2^{k-3}+1 : 2^{k-2}) = Y_{N}(1 : 2^{k-3})+2^{2};$$

And the procedure can go on until the first component of vector $C_{N}$  is reached by successive truncations like those above: 
$$C_{N}(1) = 2^{k-1}; \, Y_{N}(2) = Y_{N}(1) + 2^{k-1};$$
PROPOSITION 3.  If $N = 2^{k}$ and  $k \geq 2$, the {\it Fourier Permutation} of the indices $[1, \dots ,N]$ , has the following property:
$$\forall p \in \overline{0,(k-1)}: Y_{N}(2^{ k-p-1}+1 : 2^{k-p}) = Y_{N}(1 : 2^{k-p-1})+2^{p}$$

\vspace{5mm}

\begin{center}
\section{Direct Fourier Permutation Algorithm}
{\bf (Implicit Bit Reversal)}
\end{center}

By applying the properties mentioned in Proposition 2 and their consequences, the generating function of the the {\it Fourier Permutation} of indices $\overline{1,2^{k}}$  can be rewritten in $k$ iterations, as follows:\\

function $V=dfp(b,N)$ 

\% $N=2^{k}, k>0$;

\% $b$ is a natural number;

\% $V$ is the {\it Fourier Permutation} of indices $[ b , b+1 ,..., b+2^{k}-1 ]$

$V=[b];p2=N/2$;

while $p2 \geq 1$

\hspace{5mm} $V=[V , V+p2]$; \% one vectorial addition and one memory reallocation of $V$

\hspace{5mm} $p2=p2/2$; \% update by one division

end;\\

Let $V = dfp(0,N)$ . Taking into account that the decimal number obtained by reversing the $k$-bit representation of the decimal number $2^{k-p}$  is  $2^{(k-1)-(k-p)}=2^{p-1}$, the function that returns the decimal values corresponding to the binary representations obtained by reversing the $k$-bit representation of all of the components of $V$ is the following:\\

function $W = mkB10RevKBit(N)$ 

\% $N=2^{k}, k>0$;

$W=[0];p2=1$;

while $p2 \leq N/2$

\hspace{5mm} $W=[W , W+p2]$;

\hspace{5mm} $p2=p2*2$;

end;\\

PROPOSITION 4. On each iteration within $mkB10RevKBit$ function, the intermediate result $W$ is ascendently ordered. \\
{\it Consequence: }

The result  $W = mkB10RevKBit(N)$ is ascendently ordered. As the length of $W$ is  $N = 2^{k}$ and the values in $W$ are distinct natural numbers corresponding to binary $k$-bit representations, it follows that the maximal item in $W$ is  $2^{k}-1$. Consequently  $W = [0,1,2,3, \dots , 2^{k}-1]$. \\

The above considerations allow for the formulation of the following theorem:\\

THEOREM 1. (Correctness and Complexity of {\it Direct Fourier Permutation Algorithm}) 
If $V = dfp(b,N)$, $N = 2^{k}$, then by implication the vector $V$ meets the relation:
$$V = b+R, (Matlab formalism),$$
where $R$ is a permutation of $W$, specifically the one that corresponds to the reversed $k$-bit representations of the components of $W$ (i.e. $R=dfp(0,N)$, or in other words, reversed $k$-bit representation of $R$ is sorted in ascending order, i.e. $R=bitrevorder(W)$). The arithmetical complexity of the computation of $V$ is $O(N+k-1)$ . \\

{\it Note: } The above theorem proves that the {\it Fourier Permutation} is uniquely determined by its first component and by its most important property formulated as Proposition 3. This is because neither Proposition 3 nor Proposition 2 really depends on the first component of the index to be permuted ($[0,1,2, \dots ,2^{k-1}]$ or $[1,2, \dots ,2^{k}]$ or $[b,b+1,b+2, \dots ,b+2^{k-1}]$). In other words, all the properties mentioned in this paper regarding both the {\it Fourier Permutation} and the {\it set of additive constants of the Fourier Permutation}, including the above theorem, are independent of the particular choice of the initial index (languages like C uses 0-based indexing while Matlab uses 1-based indexing).  These are the reasons why the formalism in the above theorem has been chosen to unify the descriptions of both cases mentioned in Example 1.

\begin{center}
\section{Vectorial Digit Reversal Generalization}
\end{center}

As a natural generalization of the {\it Direct Fourier Permutation Algorithm} we propose the following algorithm that is enabled to compute the $r$-digit reversal of natural index $[1,2,3, \dots ,r^k]$ for arbitrary radices $r$ in a vectorial manner:\\

function $V=vdigitrevorder(N,r) $

\%$N=r^k$

$crN = N/r;$

$KV = crN*[0:r-1] + 1;$

$V=KV;$

$crLenV=r;$

$crStep=1;$

while $crLenV < N$

    \hspace{5mm}if $crLenV==r^{crStep}$

        \hspace{5mm}\hspace{5mm}$crStep=crStep+1;$

        \hspace{5mm}\hspace{5mm}$crN=crN/r;$

        \hspace{5mm}\hspace{5mm}$KV=V+crN;$

        \hspace{5mm}\hspace{5mm}$crLenKV=crLenV;$

        \hspace{5mm}\hspace{5mm}$crLenV=2*crLenV;$

    \hspace{5mm}else 

        \hspace{5mm}\hspace{5mm}$KV=KV+crN;$

        \hspace{5mm}\hspace{5mm}$crLenV=crLenV+crLenKV;$

    \hspace{5mm}end;

    \hspace{5mm}$V=[V,  KV];$

end;\\

The idea of the above algorithm is to compute (at each step of the iteration) the current kernel vector $KV$ and to concatenate its value at the end of the currently calculated partial result $V$ using the following rule:\\
 
While $V$ is still a partial result (i.e. length of $V<N$):
\begin{itemize}
\item if the current length of $V$ is equal to a natural power of the radix $r^p$ then 
	\begin {itemize}
		\item update the current kernel vector {\it increasing the number of its components  (up to $r^p$) and increasing the components themselves}: 

		$KV=V+r^{k-p-1}$ 

		\item update the current partial result {\it doubling the number of its components}: 
		
		$V=[V, V+r^{k-p-1}]$

	\end{itemize}
\item else (the current length of $V$ is between $r^p$ and $r^{p+1}$)
	\begin{itemize}
		\item update the current kernel vector {\it increasing all of its components}: 
	
		$KV=KV+r^{k-p-1}$ 

		\item update the current partial result {\it increasing the number of its components with $r^p$}: 
		
		$V=[V, KV]$ 

	\end{itemize}
\end{itemize}

\hspace{5mm}Another form of the same algorithm, in which the similarity to the {\it Direct Fourier Permutation Algorithm} is obvious,  is the following:

\begin{center}
\begin{tabular}{lll}

 {\bf Vectorial Digit Reversal:} &
$ $ &
 {\bf Direct Fourier Permutation:} \\

 function $V=vdro(N,r)$ &
$ $ &
 function $V=dfp(b,N)$ \\

 $V=[1]; pr=N/r;$ &
$ $ &
 $V=[b]; p2=N/2;$ \\

 while $pr \geq 1$ &
$ $ &
 while $p2 \geq 1$ \\

 \hspace{5mm}$KV=V;$ &
$ $ &
$ $ \\

  \hspace{5mm}for $cont=1:r-1$ &
$ $ &
 $ $ \\

  \hspace{10mm}$V=[V, KV+cont*pr];$ &
$ $ &
 \hspace{5mm}$V=[V , V+p2];$ \\

  \hspace{5mm}end; &
$ $ &
$ $ \\

  \hspace{5mm}$pr=pr/r;$ &
$ $ &
\hspace{5mm}$p2=p2/2;$ \\

 end; &
$ $ &
 end; \\

\end{tabular}
\end{center}

We prefer the above formulation of the {\it Vectorial Digit Reversal Algorithm} for two reasons: it is less redundant and it enables us to observe that these two functions $vdro(N,r)$, $dfp(b,N)$, have the same arithmetic complexity when $r=2$ and $N=2^k$. Also, during various tests performed by the author, the above formulation of the {\it Direct Fourier Permutation Algorithm} has been proved to be the fastest Matlab script for computing bit reversal permutation. On the other hand, both the present mathematical formulation of the bit reversal computation and the theoretical arithmetic complexity obtained in Theorem 1 suggest that the {\it Direct Fourier Permutation Algorithm} defines the minimal computational effort (minimal arithmetic complexity) for computing bit reversal permutation of an arbitrary-based index [$b, b+1, \dots, b+2^{k-1}]$ in natural arithmetic, unless a stronger property than Proposition 3 can be formulated. \\

Despite the Matlab formalism that has been chosen in the description of both of the above algorithms, up to this point of the present work, there was no particular hypothesis (concerning some particular implementation or some particular computational advantage that could have been gained by programming in one specific language, medium or platform) being pursued,  and consequently, all the above results are purely arithmetical. In the following section we will see how the particular result of Theorem 1 that concerns complexity can be refined by making the most of the speed of the vectorized Matlab calculation.

{\centering
\section{Benchmark}
}
{\centering
\subsection{Benchmark considerations}
}

In this section we aim to obtain experimentally determined time-complexity results for the two Matlab script functions $dfp$ and $vdro$, coded above. \\

{\bf Test variables:} the arithmetic complexity of the  $dfp$ function doesn't depend on the starting value ($b$) of the index [$b, b+1, \dots, b+2^{k-1}]$ and consequently the  $dfp$ function will be tested only against increasing values of the length of the index ($N=2^k$). According to the above considerations, almost identical time-complexity results are expected to be found for the $dfp(b,2^k)$ and $vdro(2^k,2)$ computations. The $vdro$ function will be tested against the variable $(k,r)$. \\

{\bf Computing the medium execution times:} the medium execution time for each test variable ($k$ and $(k,r)$ respectively)  will be the average of execution times cumulated over a great number of repetitions. Each medium execution time thus calculated will obviously depend on the performance of the computer running the tests, but the {\it nature} of variation of the medium execution time against the test variables is a characteristic of the Matlab implementations of the algorithms in themselves and does not depend on the performance of the computer. On the other hand, it is all the more necessary to establish with accuracy the medium execution time, as the execution time intended to be estimated is shorter (cases of immediate practical interest, where test variables are small). In these cases the number of repetitions will be the greatest and will decrease gradually as the value of the test variables increases, up to a value that guarantees that the resulting vector of medium execution times is a fair statistical reflection of reality. In establishing the number of repetitions corresponding to each individual test variable, as well as in establishing the minimal number of repetitions corresponding to high values of the test variables, we will consider to be a {\it fair reflection of reality} such a representation of the medium execution times in which the first symptom of convergence is present, i.e. the curve of medium execution times is nearly a smooth, ascendent one.\\

{\bf Limitations of the present {\it digitrevorder} Matlab function:} there are three reasons for replacing the existing $digitrevorder$  Matlab function within the Signal Processing Toolbox:

i) the arbitrary radix $r$ is limited to the integer range from 2 to 36. In the proposed implementation of the Vectorial Digit Reversal Algorithm ($vdro$ function) there is no such limitation. 

ii) the computation of the $vdro$ proposed function is much faster than that required to be done in the present $digitrevorder$ function.

iii) unexpected behavior of the $digitrevorder$ function could be sometimes obtained ($digitrevorder(1:81,3)$ - for example). This is because of a data validation issue: the power $k$ of the radix $r$ is computed in the digitrevorder.m (present file) as being: 
\begin{center}
$radixpow=floor(log10(N)/log10(radixbase))$ 
\end{center}
and this instruction sometimes fails to return the correct result of the calculus (see $floor(log10(27)/log10(3))$, $floor(log10(81)/log10(3))$, $floor(log10(7^7)/log10(7))$ \\
for example) and consequently, that instruction must be replaced.

\begin{center}
\subsection{Testing the proposed functions against the existing digitrevorder Matlab function}
\end{center}

This section assumes that the tested functions are called with the following syntaxes: $dfp(1,2^k)$, $vdro(2^k,2)$, $digitrevorder(1:2^k, 2)$ where $k$ takes several increasing natural values, unless other explicit syntax is present. We preferred to test $digitrevorder$ instead of $ bitrevorder$ function because, due to its  implementation, the latter is calling the former.

\begin{figure}
\centering
\includegraphics[height=9cm,width=10cm]{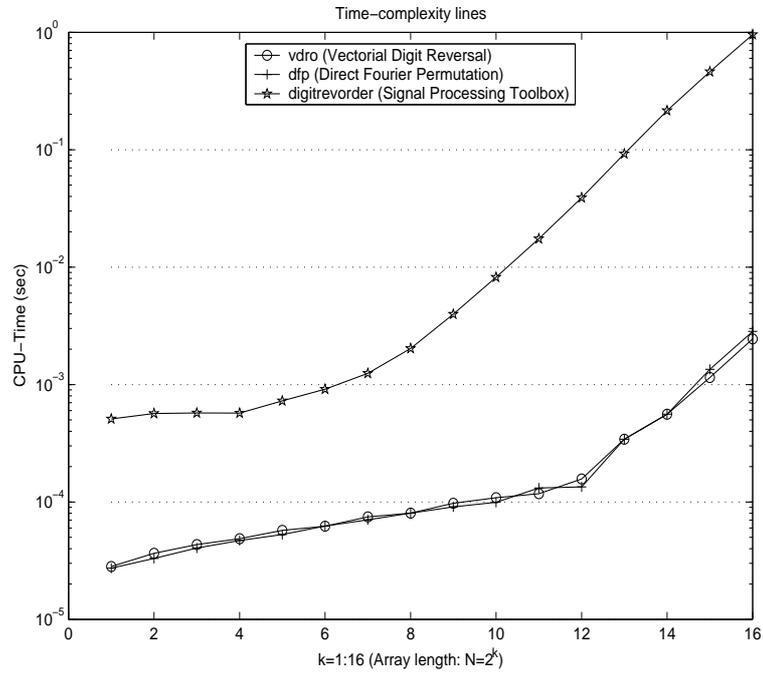}
\caption[Testing $dfp$, $vdro$, $digitrevorder$ against increasing $k$ values]
{Testing $dfp$, $vdro$, $digitrevorder$ against increasing $k$ values}
\end{figure}
\begin{figure}
\centering
\includegraphics[height=9cm,width=10cm]{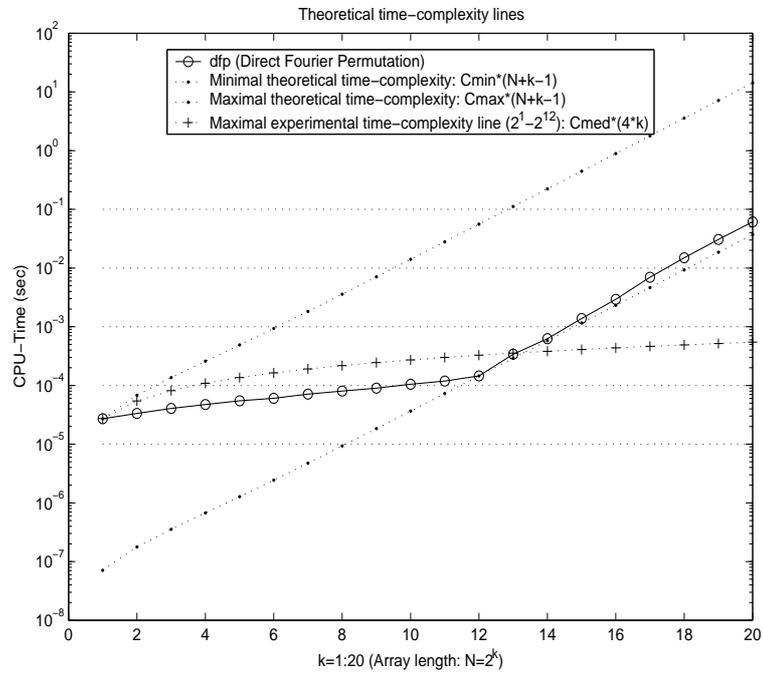}
\caption[Time-complexity lines of the $dfp$ computation]
{Time-complexity lines of the $dfp$ computation}
\end{figure}
\begin{figure}
\centering
\includegraphics[height=9cm,width=10cm]{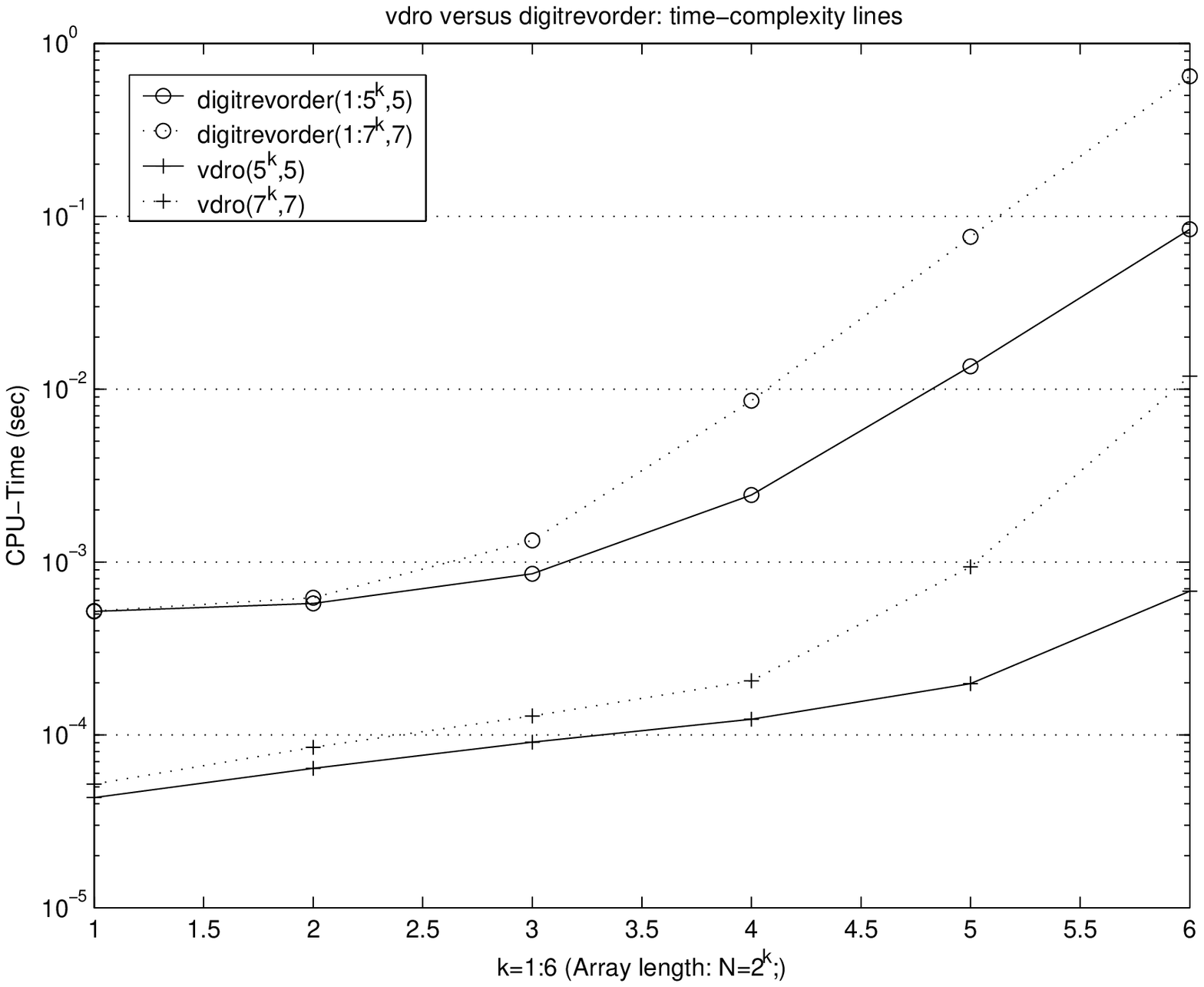}
\caption[Time-complexity lines: $vdro$ versus $digitrevorder$ (arbitrary radices)]
{Time-complexity lines: $vdro$ versus $digitrevorder$ (arbitrary radices)}
\end{figure}

As expected, the time-complexity lines (Figure 1) of the two proposed algorithms ($dfp$ and $vdro$) are almost identical when the radix is $r=2$. In this case the graphical representation in Figure 1 enables us to distinguish two variation regimes along the time-complexity lines of the both algorithms. Figure 2 reveals more clearly the point where the variation of the time-complexity line undergoes a significant change. During all the tests done by the author so far, the existence of this point has been proved to be an invariant of the algorithm. This is because the medium execution time calculated for each $k$ depends on the following factors: $k$ itself, the number of repetitions, the type of the {\it cache} and influences from other processes that are running on the same CPU. 
The greater the increase of $k$ and in the number of repetitions, the greater the contribution of the two other factors will be to the calculated medium execution time and the computation itself will becomes less cache-optimal. In Figure 2, Cmin, Cmed and Cmax are the minimal, medium, maximal duration, respectively, of the computational cycle of the $dfp$ algorithm, all of them being experimentally determined.\\

Also, Figure 2 and the $dfp$ function itself allow for the formulation of the following remarks:\\

i)	due to its simplicity, the $dfp$ function certainly has a lower complexity and a higher performance than the $bitrevorder$ function within the Signal Processing Toolbox.

ii)	the $dfp$ function allows us to compute the bit reversal order of indices  $\overline{1,2^{k}}$ in just $k$ iterations, each of these involving only three operations: one vectorial addition, one memory reallocation, and one division. 

iii)	for $k$ ranging between 2 and 12, the medium execution time needed to compute  $V = dfp(1,2^{k})$ is increasing linearily with $4k$ (reflecting the very few operations within each iteration). On this range of $k$ values, the $dfp$ function  takes all computational advantage of its simplicity. Consequently, this is the range on which the $dfp$ function reaches its maximal efficiency and its minimal computational complexity  $O(4k)$.

iv)	for $k$ ranging between 13 and 20, due to the increasing size of variable $V$, the computational complexity of the $dfp$ function suddenly turns into its own worst case scenario - predicted by its theoretical complexity $O(N+k-1)$ - on which the medium execution time needed to compute  $V = dfp(1,2^{k})$ is increasing exponentially with $k$ (i.e. linearily with  $N = 2^{k}$). But even on this range of $k$ values, the multiplicative constant characterising the time-complexity of the computation is relatively small. Consequently, on this range of $k$ values the $dfp$ function is still operative, even though it reaches its minimal efficiency and its maximal computational complexity  $O(2^{k}+k-1)$.\\

Figure 3 reveals the great improvement in the performance of the proposed $vdro$ function compared to the present Matlab function $digitrevorder$. Even for arbitrary radices $r$ the computation of the $vdro$ function is faster and more reliable.\\

{\centering
\subsection{Replicability of results}
}

 \hspace{5mm} The essential qualitative results illustrated by the tests above are the following:

i) higher performance of the $dfp$ function compared to the Matlab $bitrevorder$ function;

ii) linear variation proportional with $4k$ of the medium execution time needed to compute $V = dfp(1,2^{k})$ for $k$ ranging between $2$ and $12$.

iii) exponential variation proportional with  $(2^{k}+k-1)$ of the medium execution time needed to compute $V = dfp(1,2^{k})$  for $k$ ranging between $13$ and $20$.

iv) higher performance of the $vdro$ function compared to the Matlab $digitrevorder$ function for arbitrary radices $r$;
\\

The replication of the qualitative results mentioned above does not depend on the PC hardware architecture used. Nonetheless, the replication of any quantitative results requires a hardware architecture which must be very similar to the one used for the purposes of this paper. All the tests mentioned above were run in Matlab 6.5 R13, on the following configuration:  {\it CPU:} Prescott, Intel Pentium 4E, 2800Mhz; {\it Memory bus features:} Type - Dual DDR SDRAM, Bus width - 64 bit, Real clock - 200MHz (DDR), Effective clock - 400MHz, Bandwidth – 6400MB/s; {\it  Memory module features:} Size - 512 MB, Module type - Unbuffered, Speed - PC3200 (200MHz); \\

As far as the variation of experimental medium execution times along the range of $k$ is concerned, all the tests that have been run have high statistical relevance, being based on a number of repetitions that is high enough for statistic phenomena to become observable. All the above tests have been validated by the author by performing more tests on other hardware configurations, with similar results.\\

\newpage
\begin{center}
\section{ Appendix}
\end{center}

1.	Discrete Fourier Transform:
\\The Discrete Fourier Transform of the signal of finite length  $x(0 : N-1)$ is vector  $X(0 : N-1)$, with the following components:
$$\forall k \in \overline{0,(N-1)}: X(k) = \sum_{n=0}^{N-1}x(n)W_{N}^{kn}, \, W_{N} = exp(-2\pi i/N);$$
2.	Fast Fourier Transform, Danielson-Lanczos Lemma:
\\Fast Fourier Transform is the algorithm that allows for the computation of the Discrete Fourier Transform with a complexity of order  $O(N log_{2}N)$. The FFT algorithm for computing the Discrete Fourier Transform of signals of length  $2^{k}$ is based on the Danielson-Lanczos Lemma:
\\

{\it Danielson-Lanczos Lemma:}
\\Let signal $x(0 : N-1)$, $N = 2^{k}$ and $f(0 : N/2-1)$, $g(0 : N/2-1)$ be defined by:
$$f(n) = x(2n), \, g(n) = x(2n+1), \, n \in \overline{0,(N/2-1)};$$
Let $X(0 : N-1)$, $F(0 : N/2-1)$,  $G(0 : N/2-1)$  be the Fourier transforms of signals $x$, $f$, and $g$, respectively. 
\\Then $X$ has the following components:
\\

$X(k) = F(k) + W_{N}^{k}G(k), \, k \in \overline{0,(N/2-1)}$;

$X(N/2+k) = F(k) - W_{N}^{k}G(k), \, k \in \overline{0,(N/2-1)}$;
\\

\begin{center}
{\bf Acknowledgements}
\end{center}

The author wishes to thank professor Luminita STATE for carefully reading this text and also for all discussions and advice on Fourier Analysis and Signal Processing subjects. 

\begin{center}
{\bf REFERENCES}
\end{center}

 \char 91 1\char 93 \, Duraisamy Sundararajan, M. Omair Ahmad, M.N.S. Swamy {\it Fast Computation of the Discrete Fourier Transform of Real Data}, IEEE Transactions on Signal Processing, vol. 45, no. 8, pp. 2010-2022, Aug. 1997.

\char 91 2\char 93 \, Larry Carter and Kang Su Gatlin, {\it Towards an Optimal Bit-Reversal Permutation Program}, Proceedings of IEEE-FOCS'98, 8-11 November in Palo Alto, CA.

\char 91 3\char 93 \, Angelo A. Yong, {\it A Better FFT Bit-Reversal Algorithm Without Tables}, IEEE Transactions on Signal Processing, vol. 39, no. 10, pp. 2365-2367, Oct. 1991.

\char 91 4\char 93 \, Amitava Biswas, {\it Bit Reversal in FFT From Matrix Viewpoint}, IEEE Transactions on Signal Processing, vol. 39, no. 6, pp. 1415-1418, Jun. 1991.

\char 91 5\char 93 \, Jechang Jeong, William J. Williams, {\it A Unified Fast Recursive Algorithm for Data Shuffling in Various Orders}, IEEE Transactions on Signal Processing, vol. 40, no. 5, pp. 1091-1095, May. 1992.

\char 91 6\char 93 \, Michael Orchard, {\it Fast Bit-Reversal Algorithms Based on Index Representations in} $GF(2^{b})$, IEEE Transactions on Signal Processing, vol. 40, no. 4, pp. 1004-1007, Apr. 1992.

\char 91 7\char 93 \, Juan M. Rius, R. De Porrata-Doria, {\it New Bit-Reversal Algorithm}, IEEE Transactions on Signal Processing, vol. 43, no. 4, pp. 991-994, Apr. 1995.

\char 91 8\char 93 \, Karim Drouiche, {\it A New Efficient Computational Algorithm for Bit Reversal Mapping}, IEEE Transactions on Signal Processing, vol. 49, no. 1, pp. 251-254, Jan. 2001.

\char 91 9\char 93 \, Soo-Chang Pei, Kuo-Wei Chang, {\it Efficient Bit and Digital Reversal Algorithm Using Vector Calculation}, IEEE Transactions on Signal Processing, vol. 55, no. 3, pp. 1173-1175, Mar. 2007.

\end{document}